\documentclass[a4paper,12pt]{amsart}

\usepackage{cellspace}
\usepackage[marginparsep=0.6cm,marginparwidth=3cm]{geometry}

\usepackage{amsmath,amsthm,amsfonts,amssymb,mathtools,mathrsfs,bm}
\usepackage[utf8]{inputenc}

\usepackage[T1]{fontenc}
\usepackage{lmodern}

\usepackage{accents}
\usepackage{nicefrac}
\usepackage{diagbox}
\usepackage{makecell}

\usepackage{tikz-cd}
\usetikzlibrary{babel}
\usetikzlibrary{decorations.pathmorphing}
\usepackage[shortlabels]{enumitem}
\usepackage{float}
\usepackage{subcaption}
\usepackage{multirow,boldline}
\usepackage{graphicx}

\usepackage{xcolor}
\definecolor{colordelink}{rgb}{0,0,0.50}
\definecolor{colordecite}{rgb}{0,0.5,0}
\definecolor{colordeurl}{rgb}{0,0.41,0.5}
\usepackage[pagebackref=true]{hyperref}
\hypersetup{
    colorlinks=true,
    linkcolor=colordelink,
    citecolor=colordecite,
		urlcolor=colordeurl,
}
\usepackage[capitalize,nameinlink,noabbrev]{cleveref}
\usepackage[all]{hypcap}
\usepackage{todonotes}

\usepackage{adjustbox}

\def\im{\operatorname{im}}
\def\dim{\operatorname{dim}}

\def\sgn{\operatorname{sgn}}

\newcommand{\num}{\texttt{\#}}
\newcommand{\Alt}{\textnormal{Alt}}

\newcommand{\RR}{\mathbb{R}}

\newcommand{\CC}{\mathbb{C}}

\newcommand{\ZZ}{\mathbb{Z}}

\newcommand{\eqA}{\mathscr{A}}

\newcommand{\OO}{\mathcal{O}}

\newcommand{\Ascr}{\mathscr{A}}

\newcommand{\mfrak}{\mathfrak{m}}

\newcommand{\Ubm}{\bm{U}}
\newcommand{\Vbm}{\bm{V}}

\newcommand{\Xbm}{\bm{X}}

\newcommand{\ubm}{\bm{u}}
\newcommand{\vbm}{\bm{v}}

\newcommand{\xbm}{\bm{x}}

\newcommand{\zbm}{\bm{z}}

\makeatletter
\newcommand{\tpitchfork}{%
  \vbox{
    \baselineskip\z@skip
    \lineskip-.52ex
    \lineskiplimit\maxdimen
    \m@th
    \ialign{##\crcr\hidewidth\smash{$-$}\hidewidth\crcr$\pitchfork$\crcr}
  }%
}
\makeatother

\theoremstyle{plain}
\newtheorem{theorem}{Theorem}[section]
\newtheorem{lemma}[theorem]{Lemma}

\newtheorem*{theorem**}{Theorem\theoremnum}
\newenvironment{theorem*}[1][]{%
  \edef\theoremnum{\if\relax\detokenize{#1}\relax\else~#1\fi}
  \begin{theorem**}
}{%
  \end{theorem**}
} 
\newtheorem*{conjecture**}{Conjecture\theoremnum}
\newenvironment{conjecture*}[1][]{%
  \edef\theoremnum{\if\relax\detokenize{#1}\relax\else~#1\fi}
  \begin{conjecture**}
}{%
  \end{conjecture**}
}

\theoremstyle{definition}
\newtheorem{definition}[theorem]{Definition}
\newtheorem{proposition/definition}[theorem]{Proposition/Definition}

\newtheorem{example}[theorem]{Example}

\newtheorem{remark}[theorem]{Remark}

\Crefname{notation}{Notation}{Notations}

\theoremstyle{remark}

\begin{document}
\author{I. Breva Ribes, R. Giménez Conejero}
\title[Good real pictures of corank one map germs $(\CC^n,0)\to(\CC^{n+1},0)$]{Good real pictures of corank one map germs from the $n$-space to the $(n+1)$-space}

\address{Departament de Matemàtiques,
Universitat de Val\`encia, Campus de Burjassot, 46100 Burjassot
Spain
\newline
 Departamento de Matemáticas, Universidad Autónoma de Madrid, Ciudad Universitaria de
Cantoblanco, 28049 Madrid, Spain
\newline
Instituto de Ciencias Matemáticas, ICMAT (CSIC-UAM-UC3M-UCM), 280149, Madrid, Spain
}
\email{ignacio.breva@uv.es}
\email{Roberto.Gimenez@uv.es}

\subjclass[2020]{Primary 58K60; Secondary 14P25, 32S30} \keywords{Deformations of map germs, homology, good real perturbation, complexification, M-deformation.}

\begin{abstract}
    We study corank one  $\Ascr$-finite germs $f:(\mathbb{R}^n,0)\rightarrow (\mathbb{R}^{n+1},0)$ and their complexifications. More precisely, we study when these germs provide good real pictures of the complex germs, i.e., when there is a real deformation that has the same homology in the image (hence, homotopy) than the generic complex deformation. We give a new sufficient condition that can be computed in practice, as well as examples.
\end{abstract}

\maketitle

\section{Introduction}

This paper aims to continue the work initiated in \cite{Breva_Giménez2024} and study analytic, $\eqA$-finite  map germs $f\colon (\CC^n,0)\to(\CC^p,0)$ of corank 1 admitting good real perturbation.

Essentially, a complex map germ is said to admit a good real perturbation if the homology of the discriminant of a generic perturbation can be computed just by looking at the discriminant of certain perturbation of a real map germ, which usually there are many. As we point out in \cref{rem:homotopy}, this is equivalent to saying that the inclusion of the real into the complex is a homotopy equivalence.
The technical details of the definition are discussed in \cref{sec:grp}.






In \cite{Breva_Giménez2024} a strong necessary condition for map germs to admit good real perturbation was given:

\begin{theorem}\label{thm:anterior}
 If $f$ has corank one and admits a good real perturbation, then $\mu\big(D^k(f)\big)$ is either $0$ or $1$ provided $\dim D^k(f)>0$.
\end{theorem}

Here $D^k(f)$ denotes the set of $k$-multiple points of $f$, heuristically, the (closure of the) set in $(\CC^n)^k$ of points sharing the same image via $f$.
In the corank 1 case, these spaces are isolated complete intersection singularities (ICIS) and have a well defined Milnor number $\mu$.

The study of good real perturbations can be traced back to the work of Goryunov \cite{Goryunov1991}, which motivated that of Mond \cite{Mond1996}. 

Some classifications have been made, for instance a complete classification of map germs $(\CC^2,0)\to(\CC^3,0)$ by Marar and Mond in \cite{Marar1996}, or the fact that all map germs $(\CC,0)\to(\CC^2,0)$ admit a good real perturbation, which directly follows from the works of Guse\u{i}n-Zade \cite{GuseinZade1974} and also A'Campo \cite{ACampo1975} on finding a real perturbation of a complex plane curve exhibiting as many double points as its $\delta$-invariant. 
Other known fact is that all corank 1 map germs with $\eqA_e$-codimension 1 have good real perturbation if $n+1\geq p$ (see \cite{Cooper2002,Houston2005a}).

In this paper, we restrict to the case $p=n+1$ and compute the $(k+2)$-jet of map germs that satisfy the necessary condition given in \cref{thm:anterior}, i.e., they are candidates to have a good real perturbation. We do this in \cref{thm:K+2jet} and show which of these jets indeed have a good real perturbation in \cref{thm:necessary}.

This provides a method to determine when a map germ has good real perturbation, by just inspecting the $\eqA$-equivalence class of its $(K+2)$-jet, where $K$ is fixed, which amounts to a finite number of computations.
\newline

\cref{sec:preliminaries} is dedicated to introducing the general notation, the divided differences technique and the properties of the multiple point spaces.
\cref{sec:jet_characterization} is dedicated to the proofs of the main theorems.
\cref{sec:examples} is dedicated to showing some examples on how to apply the previous theorems to obtain map germs with good real perturbations in a given pair of dimensions.
\newline

\textbf{Acknowledgments:} R.G.C. wants to thank Ignacio Crespo Pita for continuous personal support. 

Both authors were partially supported by Grant PID2021-124577NB-I00 funded by MCIN/AEI/ 10.13039/501100011033 and by “ERDF A way of making Europe”. The first-named author was partially supported by grant UV-INV-PREDOC22-2187086, funded by Universitat de València.

\section{Preliminaries}\label{sec:preliminaries}

During most of the text, unless otherwise specified, $f\colon (\CC^n,0)\to(\CC^{n+1},0)$ will be a corank 1, $\eqA$-finite map germ, i.e.,
such that the space $T\eqA f = tf(\theta_n) + wf(\theta_{n+1})$ has finite codimension in $\theta(f)$ as $\CC$-vector space. 
Here $\theta_m$ is the space of germs of analytic vector fields in $\CC^m$ and $\theta(f)$ is the set of analytic vector fields along $f$, and $tf(\xi) = df\circ \xi$, $wf(\eta) = \eta\circ f$ for $\xi\in\theta_n$, $\eta\in\theta_{n+1}$. This, and any details regarding the theory of singular map germs, can be found in \cite{Mond2020}.

We will often use the notation $\xbm = (x_1,\ldots,x_{n-1})$.

\subsection{Multiple point spaces}

Here we present a brief explanation on the multiple point spaces of a map germ.
For a deeper view of the topic, see \cite[Chapter 9]{Mond2020}.

\begin{definition}
    Let $g\colon X\to Y$ be a finite mapping between analytic spaces, the $k$-multiple point space of $g$ is defined as:
    $$D^k(g)\coloneqq\overline{\big\{    \big(x^{(1)},\dots,x^{(k)}\big)\in X^k   :   f\left(x^{(i)}\right)=f\left(x^{(j)}\right), x^{(i)}\neq x^{(j)}  \big\}}.$$
    If $f\colon(\CC^n,0)\to(\CC^{p},0)$ is a smooth, $\eqA$-finite map germ, then it admits a stable unfolding $F(x,u) = (f_u(x),u)$ with $f_0\equiv f$.
    Then, $D^k(f)$ is the analytic space defined as the germ in $(\CC^n)^k$ of
    $$D^k(f) = D^k(F) \cap \lbrace u= 0\rbrace.$$
\end{definition}

\begin{remark}\label{rem:divided_differences}
In the case that $f$ has corank $1$, each $D^k(f)$ is actually an ICIS, and one can compute their explicit equations via the method of the divided differences.
We show this method for our case $p=n+1$ and, since it is of corank 1, we can consider that $f$ is of the form
$$f(\xbm,z) = (\xbm,f_1(\xbm,z),f_2(\xbm,z))$$
writing $\xbm = (x_1,\ldots,x_{n-1})\in\CC^{n-1}$, $z\in\CC$,
for some $f_1,f_2\colon(\CC^n,0)\to(\CC,0)$.

Set $f_j^1 = f_j$ for $j=1,2$ and for each  $i\geq 1$ define
$$f_j^{i+1}(\xbm,z_1,\ldots,z_k) =\frac{f_j^i(\xbm,z_1,\ldots,z_{k-2},z_{k-1})-f_j^i(\xbm,z_1,\ldots,z_{k-2},z_k)}{z_{k-1}-z_k}.$$
These are all well-defined functions and every $D^{k}(f)$ is isomorphic to the germ of the zero set in $\CC^{n+k-1}$ cut out by $f_1^2,\ldots,f_1^k,f_2^2,\ldots,f_2^k$, for each $k\geq 2$.
Moreover each $f_j^{i+1}$ is invariant under any permutation of the $z_1,\ldots,z_{k}$.
For more on the divided differences technique we refer to \cite[end of p. 371]{Mond1987} or \cite[p. 52]{GimenezConejero2022}.

In practice, it will be very convenient to adopt the notation $S_j(z_1,\ldots,z_k)$ to describe the polynomial defined as the sum of all monomials of degree $j$ in $z_1,\ldots,z_k$.
Since we will need to use $S_j$ for different numbers of variables (deppending on the $D^k(f)$ we are working with), we will often use the notation $S_j(\zbm)$, and the length of $\zbm$ will be clear from the context.
This abuse of notation will not interfere with any of the computations.
\end{remark}

\begin{example}
    Let us show an example on how to compute equations for $D^k(f)$.
    Let $f\colon(\CC^2,0)\to(\CC^3,0)$ be the map germ given by $f(x,z) = (x,z^2,z^{2k+1}+x^mz)$ for $k,m\geq 1$.

    Here we have $f_1(x,z) = z^2, f_2(x,z) = z^{2k+1}+x^mz$.
    In this case, $D^2(f)$ is isomorphic to the ICIS in $\CC^{3}$ cut out by the equations
    \begin{align*}
        0=f_1^2(x,z_1,z_2) &= \frac{z_1^2- z_2^2}{z_1- z_2} = z_1+ z_2,\\ 
        0=f_2^2(x,z_1,z_2) &= \frac{z_1^{2k+1}- z_2^{2k+1}+x^m(z_1 -z_2)}{z_1-z_2} = \sum_{i=0}^{2k} z_1^{2k-i}z_2^i +x^m.
    \end{align*}
    In particular, it is isomorphic to the ICIS in $\CC^2$ cut out by the equation $z_1^k + x^m = 0$.
    Moreover $D^3(f)$ is empty for each $k,m\geq 1$, since one of the equations of $D^3(f)$ is
    $$0=f_1^3(x,z_1,z_2,z_3) = \frac{z_1+z_2 - (z_1 + z_3)}{z_2 - z_3} = 1$$
    and therefore every $D^k(f)$ is empty for $k\geq 3$.
\end{example}
\medskip

Multiple point spaces have evident symmetry by permutation of the components of each of their points. 
This translates into several nice properties, the first one we are interested in is the fact that $\eqA$-finiteness of $f$ can be detected by studying all of the $D^k(f)$.

Let $\Sigma_k$ be the group of permutations of order $k$.
Each permutation $\sigma$ has a unique decomposition into cycles, and the lengths of the cycles determine a partition of $k$.
For instance, if $k=7$ the permutation $\sigma = (1\ 2\ 4\ 3) (5 \ 7) (6) $ is associated to the partition $(4,2,1)$ of $7$.
Let $\sigma^\num$  be the number of cycles of the unique decomposition of $\sigma$ (in this example, $\sigma^\num = 3$).

\begin{lemma}[See {\cite[Corollary 2.15]{Marar1989}} and {\cite[Corollary 2.8]{Houston2010}}]\label{lem: dkgamma}
If $f:(\CC^n,S)\rightarrow (\CC^p,0)$ is a finite germ of corank $1$, $n < p$, and $\sigma\in\Sigma_k$, then:
\begin{enumerate}
	\item If $f$ is stable, $D^k(f)^\sigma$ is smooth of dimension $p-k(p-n)-k+\sigma^\num$, or empty.\label{iPD}
	\item $\eqA_e-codim(f)$ is finite if, and only if:\label{iiPD}
	\begin{enumerate}
		\item for each $k$ with $p-k(p-n)-k+\sigma^\num\geq0$, $D^k(f)^\sigma$ is empty or an ICIS of dimension $p-k(p-n)-k+\sigma^\num$,
		\item for each $k$ with $p-k(p-n)-k+\sigma^\num<0$, $D^k(f)^\sigma$ is a subset of $S^k$, possibly empty.
	\end{enumerate}
\end{enumerate}
We say that $d_k\coloneqq p-k(p-n)$ is the {\it expected dimension} of $D^k(f)$ and that $d_k^\sigma\coloneqq p-k(p-n)-k+\sigma^\num$ is the {\it expected dimension} of $D^k(f)^\sigma$.
\label{DP}
\end{lemma}

\begin{example}
    In our previous example, for $f(x,z) = (x,z^2,z^{2k+1}+x^mz)$ for $k,m\geq 1$, we computed that $D^2(f)$ is isomorphic to the zero-set of $z^k+x^m$ in $\CC^2$.
   Hence, $D^2(f)$ is nonempty of dimension 1 for all $k,m\geq 1$, which is the expected dimension, and $D^k(f)$ is empty for $k\geq 3$.
   By the previous theorem, $f$ is $\eqA$-finite and is moreover stable if and only if $m=1$ since $D^2(f)$ is smooth in that case.
\end{example}

The following result summarizes the relation between singularity and dimensions of the multiple point spaces, it is \cite[Theorem 4.8]{GimenezConejero2022c}.

\begin{theorem}\label{thm:dk_singular}
    If $f\colon(\CC^n,S)\to(\CC^p,0)$ is an unstable, $\eqA$-finite map germ of corank one, $n<p$ and spaces $D^k(f)\neq \emptyset$ with $d_k\geq 0$, then the following are equivalent:
    \begin{enumerate}
        \item $D^k(f)$ has a singularity;
        \item $D^k(f_s)$ has non-trivial homology in the middle dimension, $f_s$ stable;
        \item $D^{k+j}(f)$ has a singularity if it is not empty, for every $j\geq 0$ such that $d_{k+j}\geq 0$;
        \item every space $D^k(f)^\sigma$ has a singularity if it is not empty and $d_k^\sigma \geq 0$;
        \item some space $D^k(f)^\sigma$ has a singularity if $d_k^\sigma\geq 0$;
        \item $D^k(f_s)$ has non-trivial alternating homology in middle dimension.
    \end{enumerate}
    Furthermore, if $d_k^\sigma < 0$ and $D^k(f)^\sigma$ is not empty, all of the above are satisfied.
\end{theorem}

The following is a well-known and useful lemma.
A proof can be found, for instance, in \cite[Lemma 2.2]{Breva_Giménez2024}

\begin{lemma}\label{lem:ICIS_A1}
   Let $(X,0)$ be an isolated complete intersection singularity (ICIS) of positive dimension.
   If $\mu(X,0) = 1$, then $(X,0)$ is isomorphic to an $A_1$ singularity, i. e., its local algebra is isomorphic to the local algebra of $V(x_1^2+\cdots+x_n^2)$ in $(\CC^n,0)$.
\end{lemma}

Another interesting property of multiple point spaces is that they can be used to compute the homology in the image of a finite map.
To show this result we need to introduce alternating homology, which also uses the symmetry of $D^k(f)$.

\begin{definition}
    Define
    \begin{align*}  
   C^\Alt_n(D^k(f); \ZZ) &\coloneqq \lbrace c\in C_n(D^k(f);\ZZ)\ |\ \sigma(c) = \sgn(\sigma)c \, \forall \sigma\in\Sigma_k\rbrace, \\
   C^\Alt_n(D^k(f); G) & \coloneqq C^\Alt_n(D^k(f); \ZZ)\otimes G
    \end{align*}
    for $G$ an abelian group.
    One can check that these determine a chain complex with boundary defined by the restriction of the usual chains.
    
    The alternating homology groups of $D^k(f)$, denoted $AH_*(D^k(f))$, are then defined as the homology groups of $C^\Alt_*(D^k(f))$ with respect to the border map restriction.
\end{definition}

We then have the following result, of which several versions exist for different contexts, we use the one in \cite[Theorem 5.4]{Houston2007}, which is very general (see also \cite{Goryunov1993,Goryunov1995,CisnerosMolina2022,Mond2020}). 

\begin{theorem}\label{thm:icss}
    Let $g\colon X\to Y$ be a finite, surjective and analytic map between compact subanalytic spaces. 
    Let $G$ be an abelian group.
    
    The spectral sequence whose first page is defined as $E_{p,q}^1(g;G) = AH_q(D^{p+1}(g); G)$ and with differentials induced by the natural map $\pi\colon D^{k+1}(f)\to D^k(f)$ converges:
    $$E_{p,q}^1(g; G) \implies H_*(g(X); G).$$
\end{theorem}

\subsection{Good real perturbations}\label{sec:grp}

In this section we recall the concept of good real perturbation and some related results.

\begin{definition}
    A map $g\colon X \to Y$ is {\it locally stable} if its restriction to its set of non-submersive points is finite and $g\colon(X,g^{-1}(y))\to(Y,y)$ is stable for every $y\in Y$.
    
    A {\it stabilisation} of $f\colon(\CC^n,0)\to(\CC^{n+1},0)$ is a 1-parameter stable unfolding $F(x,t) = (f_t(x),t)$ admitting a representative such that $f_t$ is locally stable for all $t$.
    We will also say that $f_t$ is a stabilisation of $f$.
\end{definition}

\begin{definition}
    Let $f\colon(\CC^n,0)\to(\CC^{n+1},0)$ and $f^\RR\colon(\RR^n,0)\to(\RR^{n+1},0)$ be smooth map germs such that the complexification of $f^\RR$ is $\eqA$-equivalent to $f$.
    Let $f_t, f_t^\RR$ be stabilisations of their respective germs.
    We say that $f_t^\RR$ is a good real perturbation of $f$ if
    $$\beta_i\left(\im f_t\right) \neq 0 \implies \beta_i\left(\im f_t\right)=\beta_i\left(\im f_t^\RR\right).$$
    
\end{definition}

\begin{remark}\label{rem:homotopy}
In \cite[Corollary 3.7]{Breva_Giménez2024} it is proven that if $f^\RR_t$ is a good real perturbation of $f$, then $\beta_i\left(\im f_t\right) = \beta_i\left(\im f^\RR_t\right)$. More generally, for the kind of germs we study it is proven that if $f_t^\RR$ is a good real perturbation then the inclusion of the real image into the complex is a homotopy equivalence, see \cite[Theorem 5.4]{Breva_Giménez2024}.
\end{remark}

Some classifications of map germs with good real perturbations have been carried out, for instance the case of curves, $(\CC,0)\to(\CC^2,0)$, was solved by Guse\u{i}n-Zade \cite{GuseinZade1974} and A'Campo \cite{ACampo1975}, who proved that every complex plane curve admits a real perturbation presenting as many real nodes as its $\delta$-invariant.

In \cite{Marar1996} the classification of map germs that admit good real perturbation in $(\CC^2,0)\to(\CC^3,0)$ was carried out, showing that the only map germs with good real perturbation in these dimensions are
$$(x,y^3,xy+y^{3k-1}), \, k\geq 1$$
which coincide with $S_1$ for $k=1$ and with $H_k$ for $k\geq 2$.

In \cite[Theorem 4.14]{Breva_Giménez2024}, a necessary condition for a map germ $f\colon(\CC^n,0)\to(\CC^{n+1},0)$ to have good real perturbation is given, we stated it in the introduction but we recall it here for the sake of completion:

\begin{theorem}\label{thm:grp_necessary}
If $f\colon(\CC^n,0)\to(\CC^{n+1},0)$ is of corank 1 and admits a good real perturbation, then $\mu(D^k(f)) \leq 1$ whenever $\dim D^k(f) \neq 0$.
\end{theorem}

\begin{remark}\label{rem:althominreal}
Many lemmas were involved in the proof of this theorem and one that is very useful for us is that the alternating homology of the multiple point spaces of the real and the complex map must be equal to have a good real picture (\cite[Theorem 4.3]{Breva_Giménez2024}). This is, indeed, also sufficient by \cref{thm:icss}.
\end{remark}

The next result states that for $p=n+1$ and $n\geq 3$, candidates for good real perturbations cannot present non-empty multiple point spaces of dimension zero, which could only be $D^{n+1}(f)$.

\begin{theorem}\label{thm:dimension}
If $f\colon(\CC^n,0)\to(\CC^{n+1},0)$ is $\eqA$-finite of corank 1, and $\mu(D^k(f)) \leq 1$ for those multiple point spaces of positive dimension, then $D^{k}(f) = \emptyset$ for every $k > \lceil \frac{n+3}{2} \rceil$.

In particular, if $n\geq 3$, then the multiple point space with expected dimension $0$, $D^{n+1}(f)$, must be empty.
\end{theorem}
\begin{proof}

Since $f$ is of corank 1, $D^k(f)$ is isomorphic to an ICIS in $(\CC^{n+k-1},0)$ whose equations can be obtained via divided differences as we explained in \cref{rem:divided_differences}, from which we will be using the notation from \cref{rem:divided_differences} (except for the name of the variables), that is, $f_1^2,\ldots,f_1^k,f_2^2,\ldots,f_2^k$ are the equations in $\CC^{n+k-1}$ of an ICIS isomorphic to $D^k(f)$, where $f_1^i$, and $f_2^i$ depend on variables $x_1,\ldots,x_n,x_{n+1},\ldots,x_{n+i-1}$ and are invariant under any permutation of the variables $x_n,\ldots,x_{n+i-1}$.

To prove our claim, assume that for $k_0:= \lceil \frac{n+3}{2}\rceil$ we have $D^{k_0}(f)\neq\emptyset$ (otherwise, the result is obvious).

If $n=2$, then $k_0=3$ and necessarily $D^k(f) = \varnothing$ for $k>3$ since $f$ is $\eqA$-finite.

Assume $n>3$.
Since $f$ is $\eqA$-finite $D^{k_0}(f)$ must have positive dimension and must also have either Milnor number $1$ or $0$.
If $\mu(D^{k_0}(f)) = 0$, all $2k_0-2$ equations of $D^{k_0}(f)$ determine a smooth hypersurface, hence have some linear term. 
In particular, for each $i=2,\ldots,k_0$, $j=1,2$, we can express each equation $f^i_j = 0$ as
$$f^i_j = \alpha_{i,j}x_{s_{i,j}} + \tilde f^i_j = 0$$
for some non-zero $\alpha_i^j\in\CC$ and some $\tilde f^i_j \in\OO_{n+k_0-1}$ with no linear part in $x_{s_{i,j}}$, $s_{i,j}\in\lbrace 1,\ldots,n+k_0-1\rbrace$ with $s_{i,j}\neq s_{a,b}$ for $(i,j)\neq(a,b)$. This is because, if $s_{i,j} = s_{a,b}$, then we can substitute either $f^i_j=0$ or $f^a_b = 0$ by the equation $\tilde f_i^j - \frac{\alpha_{i,j}}{\alpha_{a,b}} \tilde f_a^b =0 $, which cannot be a trivial equation since $D^{k_0}(f)$ is an ICIS, and therefore must again be a hypersurface with a linear variable, different from $x_{s_{i,j}}$.

Since $2k_0 -2 \geq n$ and $D^{k_0}(f)$ is an ICIS in $\CC^{n+k_0-1}$, at least one of the equations must have  a term given by one of the $x_n,\ldots,x_{n+k_0-1}$ linearly.
Let $i_0, j_0$ be such that $f_{j_0}^{i_0}$ is that equation, which is invariant by permutations of the $x_n,\ldots,x_{n+k_0-1}$ variables, so we can assume it is of the form
$$f_{j_0}^{i_0} = \alpha_{i_0,j_0}(x_n+\cdots+x_{n+k_0-1}) + \tilde f_{j_0}^{i_0}$$
with $\alpha_{i_0,j_0}\neq 0$, and $\tilde f_{j_0}^{i_0}$ without linear terms on $x_n,\ldots,x_{n+k_0-1}$.  Now, the divided differences of this function give an equation
$$\alpha_{i_0,j_0} + \bar f_{j_0}^{i_0+1} = 0$$
for some $\bar f_{j_0}^{i_0+1} \in\OO_{n+k_0}$.
Moreover, $\bar f_{j_0}^{i_0+1}(0) = 0$ because $\bar f_{j_0}^{i_0}$ has no linear term in any of the $x_n,\ldots,x_{n+k_0-1}$, hence there are no solutions.
But since we assumed that $D^{k_0}(f)\neq \varnothing$, necessarily $i_0 = k_0$ and $D^{k_0+1}(f) = \varnothing$ as required.

If $\mu(D^{k_0}(f)) = 1$, then $2k_0-3$ of the $2k_0-2$ equations define a hypersurface, and the other one defines a Morse singularity. Since $2k_0 -3\geq n$, the same argument shows the result.
\end{proof}

\section{Jet characterization}\label{sec:jet_characterization}

Let $f:(\CC^n,0)\to(\CC^{n+1},0)$ be a corank one map germ that is $\Ascr$-finite. Since it has corank one, it is $\Ascr$-equivalent to a germ (we keep denoting it as $f$)
\begin{equation}\label{eq:fcorank1}
f(\xbm,z)=\big(\xbm,f_1(\xbm,z),f_2(\xbm,z)\big).\end{equation}

Moreover, if the map germ admits a good real perturbation, we know by \cref{thm:dimension} together with \cref{thm:grp_necessary,lem:ICIS_A1} that the multiple point spaces $D^k(f)$ of $f$ must be of positive dimension and either smooth, an $A_1$ singularity, or empty. As we have discussed in \cref{lem:ICIS_A1}, since $A_1$ singularities are equivalent (that is, its local algebra) to a hypersurface singularity, all the defining equations of $D^k(f)$ are smooth, except one.

 We show now that this points in the direction of a sufficient condition to have a good real perturbation. The technical conditions we impose as hypothesis will become clear when we prove in \cref{thm:necessary} that this is, indeed, a characterization of having good real perturbations (in corank one) and those hypotheses come naturally.

\begin{theorem}\label{thm:K+2jet}
Let $g:(\CC^n,0)\to(\CC^{n+1},0)$ be a germ $\Ascr$-equivalent to a real germ $f$ whose $(K+2)$-jet $j^{K+2}(f)$ is one of the following two cases:
\begin{enumerate}
\item for $\ubm\in\CC^{K-1}$, $\vbm\in\CC^{K-1}$, $\xbm\in\CC^{n-2K+1}$ and $\alpha<1$ (possibly zero)
    \[\left(\ubm;\:  \vbm;\:  \xbm; \sum_{i=1}^{K-1}z^iu_i+z^Kq(\ubm,\vbm,\xbm)+\alpha z^K v_{K-1}+z^{K+2};\ \sum_{i=1}^{K-1}z^iv_i+z^{K+1}\right);\]
    \item for $\ubm\in\CC^{K-1}$, $\vbm\in\CC^{K}$ and $\xbm\in\CC^{n-2K}$
    \[\left(\ubm;\:  \vbm;\:  \xbm; \sum_{i=1}^{K-1}z^iu_i+z^Kq(\ubm,\vbm,\xbm)+z^{K+2};\ \sum_{i=1}^{K}z^iv_i+z^{K+1}x_1\right);\]
\end{enumerate}
such that $q$ does not have a linear term in $v_{K-1}$ in the first case or $v_K$ in the second, and $q(\bm{0},\bm{0},\xbm)=\sum_i x_i^2$
or simply zero if there are no variables $\xbm$. Then, $g$ admits a good real perturbation and, indeed, for a $0<t\ll1$ it is the perturbation
\[f+(\bm{0};\: \bm{0};\: \bm{0};\: tz^K;\: 0).\]
\end{theorem}
\begin{proof}
A set of defining equations of the multiple point spaces $D^k(f)$ are the corresponding pairs of divided differences as we mentioned in \cref{rem:divided_differences}. 
In particular, if $k< K$, all the pairs of equations have the form 
\begin{equation}\label{eq:eqsK-1}\begin{aligned}
&\begin{cases}
u_i+\sum_{j=1}^{K-1-i} S_{j}(\zbm)u_{i+j}+S_{K-i}(\zbm)\big(q(\ubm,\vbm,\xbm)+\alpha v_{K-1}\big)+S_{K+2-i}(\zbm)&=0\\
v_i+\sum_{j=1}^{K-1-i} S_{j}(\zbm)v_{i+j}+S_{K+1-i}(\zbm)&=0;
\end{cases}\\
&\begin{cases}
u_i+\sum_{j=1}^{K-1-i} S_{j}(\zbm)u_{i+j}+S_{K-i}(\zbm)q(\ubm,\vbm,\xbm)+S_{K+2-i}(\zbm)&=0\\
v_i+\sum_{j=1}^{K-i} S_{j}(\zbm)v_{i+j}+S_{K+1-i}(\zbm)x_1&=0;
\end{cases}
\end{aligned}\end{equation}
where the first bracket would correspond to the equations in the first case of the statement and the second bracket to the second case, $i=1,\dots,k-1$ and each of the $S_j(\zbm)$ is the homogeneous polynomial given by the sum of all the monomials of degree $j$ in the variables $\zbm=(z_1,\dots,z_{i})$ with coefficient one (recall the abuse of notation from \cref{rem:divided_differences}). 

For $D^K(f)$ we add the following pair equations to all the previous defining equations (again, depending on the case):
\begin{equation}\label{eq:eqsK}\begin{aligned}
&\begin{cases}
u_{K-1}+(z_1+\cdots+z_{K-1})\big(q(\ubm,\vbm,\xbm)+\alpha v_{K-1}\big)+S_{3}(z_1,\dots,z_{K-1})&=0\\
v_{K-1}+S_{2}(z_1,\dots,z_{K-1})&=0;
\end{cases}\\
&\begin{cases}
u_{K-1}+(z_1+\cdots+z_{K-1})q(\ubm,\vbm,\xbm)+S_{3}(z_1,\dots,z_{K-1})&=0\\
v_{K-1}+(z_1+\cdots+z_{K-1})v_K+S_2(z_1,\dots,z_{K-1})x_1&=0.
\end{cases}
\end{aligned}\end{equation}

Similarly, to obtain the defining equations of $D^{K+1}(f)$, we take the equations of $D^K(f)$ plus a new pair of equations for each case, given by
\begin{equation}\label{eq:eqsK+1}\begin{aligned}
&\begin{cases}
q(\ubm,\vbm,\xbm)+\alpha v_{K-1} + S_{2}(z_1,\dots,z_K)&=0\\
z_1+\cdots+z_K&=0;
\end{cases}\\
&\begin{cases}
q(\ubm,\vbm,\xbm)+S_{2}(z_1,\dots,z_K)&=0\\
v_K+x_1(z_1+\cdots+z_K)&=0.
\end{cases}
\end{aligned}\end{equation}

Finally, for $D^{K+2}(f)$, we (would) have the pair of equations
\begin{equation}\label{eq:eqsK+2}\begin{aligned}
&\begin{cases}
z_1+\cdots+z_{K+1}&=0\\
1&=0;
\end{cases}\\
&\begin{cases}
z_1+\cdots+z_{K+1}&=0\\
x_1&=0;
\end{cases}
\end{aligned}\end{equation}
from which we see that $D^{K+2}(f)$ is empty in the first case. Clearly, $D^{K+3}(f)$ is empty in both cases.
\newline

We show now that all the spaces $D^k(f)$ are smooth, an $A_1$ singularity or empty. 
Indeed, we see from \cref{eq:eqsK-1,eq:eqsK} that all the spaces $D^k(f)$ with $k\leq K$ are smooth, in both cases. If $q$ contains a constant then $D^{K+1}$ is empty, in both cases, and the map germ is stable by \cref{lem: dkgamma}. Otherwise, using \cref{eq:eqsK} we can eliminate all the variables $\ubm$ and $\vbm$ by terms of order greater than two, excluding $v_{K-1}$ in the first case and $v_K$ in the second, which are substituted by terms of order two. After the corresponding substitution and elimination of the variables $\ubm$ and $\vbm$ in the equations from \cref{eq:eqsK+1},the remaining singular equation,
\begin{equation} \label{eq:q_eq_cases}
\begin{aligned}
q(\ubm,\vbm,\xbm)+\alpha v_{K-1} + S_{2}(z_1,\dots,z_K)&=0\text{ for the first case, or}\\
q(\ubm,\vbm,\xbm)+S_{2}(z_1,\dots,z_K)&=0 \text{ for the second},
\end{aligned}
\end{equation}
 will yield a Morse singularity of index zero. Clearly, if it is a Morse singularity, all the variables $\ubm$ and $\vbm$ (except $v_{K-1}$ or $v_K$) play no role since they will contribute to the singular equation with terms of order three or higher. Regarding $v_{K-1}$ or $v_K$, they do not appear linearly in $q$ in their respective cases by assumption but, in the first case, we also have the relations from \cref{eq:eqsK,eq:eqsK+1}:
\[\begin{aligned}
    v_{K-1}&=-S_2(z_1,\dots,z_{K-1})\text{, and}\\
    z_K&=-z_1-\cdots-z_{K-1};
\end{aligned}\]
and on the second case we have the relation from \cref{eq:eqsK+1}:
\[v_K=-x_1(z_1+\cdots+z_K).\]
This last one does not modify the properties of $q$ since it does not have linear term on $v_k$ by assumption.
 These yield (up to degree two), respectively for the first and second case of \cref{eq:q_eq_cases}: 
\[\begin{aligned}
q(\bm{0},\bm{0},\xbm)-\alpha S_2(z_1,\dots,z_{K-1}) + S_{2}(z_1,\dots,z_{K-1},-z_1-\cdots-z_{K-1})=&0\\
q(\bm{0},\bm{0},\xbm)+ S_2(z_1,\dots,z_{K})=&0.
\end{aligned}\]
Observe now that
\[\begin{aligned}
    S_{2}(z_1,\dots,z_{K-1},-z_1-\cdots-z_{K-1})=&S_2(z_1,\dots,z_{K-1})-\sum_{i=1}^{K-1}z_i(z_1+\dots+z_{K-1})\\&+(-z_1-\cdots-z_{K-1})^2\\
=&S_2(z_1,\dots,z_{K-1}).
\end{aligned}\]
Hence, the equations in \cref{eq:q_eq_cases} are respectively for each case, up to degree 2:
\[\begin{aligned}q(\bm{0},\bm{0},\xbm)+(1-\alpha) S_2(z_1,\dots,z_{K-1})=&0\text{, and}\\
q(\bm{0},\bm{0},\xbm)+ S_2(z_1,\dots,z_{K})=&0.\end{aligned}\]
Notice that the number of variables is not the same in both cases, but neither is the ambient space nor the dimensions. Both equations define a Morse singularity of index zero since $q(\bm{0},\bm{0},\xbm)$ is a positive definite quadratic form, $\alpha<1$ and we can perform a change of variables in the $z_i$ that makes $S_2$ equivalent to a sum of squares. 

Finally, the space $D^{K+2}(f)$ is empty in the first case and in the second case we have to add the equations
\[\begin{cases}
z_1+\cdots+z_{K+1}&=0\\
x_1&=0
\end{cases}\]
from \cref{eq:eqsK+2} to
\[q(\bm{0},\bm{0},\xbm)+ S_2(z_1,\dots,z_{K-1})=0,\]
which simply eliminates $x_1$ and $z_{K+1}$ and yields another Morse singularity of index zero.
\newline

To conclude the proof, we have to show that the deformation
\[f+t(\bm{0};\: \bm{0};\: \bm{0};\: z^K;\: 0)\]
yields the correct deformation on the multiple point spaces that makes it a good real deformation. This is easy, since that deformation has $D^{K+1}(f)$ defined by
\[\begin{aligned}
q(\ubm,\vbm,\xbm)+\alpha v_{K-1} + S_{2}(z_1,\dots,z_K)&=t\text{, or}\\
q(\ubm,\vbm,\xbm)+S_{2}(z_1,\dots,z_K)&=t.
\end{aligned}\]
We can do a similar substitution of variables as we did before using the equations from \cref{eq:eqsK}, now modified by the addition of some new terms depending on $t$. 
Observe that the  new terms of degree one or two that appear multiplied by $t$ in the previous equations (e.g., possibly $u_{K-1}$ if it appears linearly in $q$) will not affect the index of the remaining quadratic form, since $t$ is arbitrarily small. 
Hence, we will have that $D^{K+1}(f)$ defines a variety isomorphic to
\[\begin{aligned}q(\bm{0},\bm{0},\xbm)+(1-\alpha) S_2(z_1,\dots,z_{K-1})=&t\text{, or}\\
q(\bm{0},\bm{0},\xbm)+ S_2(z_1,\dots,z_{K})=&t.\end{aligned}\]
This is diffeomorphic to a sphere and, in addition, the permutation group acts by reflection on hyperplanes, so the top-dimensional homology is alternating (a reflection of a hyperplane changes the orientation of the fundamental class). The same situation happens for $D^{K+2}(f)$. Now, it is clear that the perturbation of the statement provides the good real perturbation (recall \cref{rem:althominreal}).
\end{proof}

\begin{theorem}\label{thm:necessary}
Any $\Ascr$-finite corank one map germ $g:(\CC^n,0)\to(\CC^{n+1})$, $n>2$, that has a good real picture is equivalent to a real map germ $f$ whose $(K+2)$-jet $j^{K+2}(f)$ is one of the cases in \cref{thm:K+2jet}, i.e., one of the following two cases:
\begin{enumerate}
    \item for $\ubm\in\CC^{K-1}$, $\vbm\in\CC^{K-1}$, $\xbm\in\CC^{n-2K+1}$ and $\alpha<1$ (possibly zero)
    \[\left(\ubm;\:  \vbm;\:  \xbm; \sum_{i=1}^{K-1}z^iu_i+z^Kq(\ubm,\vbm,\xbm)+\alpha z^K v_{K-1}+z^{K+2};\ \sum_{i=1}^{K-1}z^iv_i+z^{K+1}\right);\]
    \item for $\ubm\in\CC^{K-1}$, $\vbm\in\CC^{K}$ and $\xbm\in\CC^{n-2K}$
    \[\left(\ubm;\:  \vbm;\:  \xbm; \sum_{i=1}^{K-1}z^iu_i+z^Kq(\ubm,\vbm,\xbm)+z^{K+2};\ \sum_{i=1}^{K}z^iv_i+z^{K+1}x_1\right);\]
\end{enumerate}
such that $q$ does not have a linear term in $v_{K-1}$ in the first case or $v_K$ in the second, and $q(\bm{0},\bm{0},\xbm)=\sum_i x_i^2$
or simply zero if there are no variables $\xbm$
\end{theorem}
\begin{proof}
Since $n>2$ and $g$ admits a good real picture, we know from \cref{thm:grp_necessary,thm:dimension} that its nonempty multiple point spaces $D^k(g)$ have positive dimension and have Milnor number zero or one. 

During all of this proof we will again make use of the explicit description of the equations describing $D^k(g)$ given by the method of the divided differences.
More specifically, we will rely heavily on the fact that the actual map germ can be easily recovered just from these equations: it suffices to look at the coefficient of each symmetric polynomial $S_j(\zbm)$ in the $\zbm$ variables, to obtain the coefficient of the corresponding $z^{k+j}$.

Let us assume that $K+1$ is the minimum multiplicity $k$ of those multiple point spaces $D^k(g)$ with Milnor number one. 
Then, $D^K(g)$ is the last multiple point space that is smooth and, hence, the defining set of equations given by the divided differences (see \cref{rem:divided_differences}) have linear part in every equation.
Moreover, the coefficients of these linear terms define a linear, invertible map in $2(K-1)$ variables. 
After a (simultaneous) change of coordinates in both the map germ $g$ and the equations of $D^K(g)$ we can assume that only one of those variables appear linearly in each equation, let us call those variables $\ubm=(u_1,\dots,u_{K-1})$ and $v_1,\dots,v_{K-1}$. 
We will rename one of the $\xbm$ to $v_K$ in one of the cases later on.

Denoting momentarily the remaining variables as $\xbm$ and $z$, we have just shown that $j^{K}g$ is $\Ascr^k$-equivalent to  
\[\left(\ubm;\  \vbm;\  \xbm; \sum_{i=1}^{K-1}z^iu_i;\ \sum_{i=1}^{K-1}z^iv_i\right).\]
Now, we know that $D^{K+1}(g)$ has Milnor number one, so it is an $A_1$ singularity by \cref{lem:ICIS_A1}. In particular, to the previous defining equations, we have to add one smooth equation and another one that gives a nondegenerate Morse singularity. 
We write them now in the most general expression as they could appear with the method of the divided differences, and we will later find equivalences that simplify these expressions. Let us now distinguish between two cases depending on what the smooth equation is, assuming as well that the smooth one is the second one:
\begin{enumerate}
    \item If the smooth equation is linear only on the $z_i$ and, maybe, on some of $\ubm$ and $\vbm$, we have an equation 
    \begin{multline*}(z_1+\dots+z_K)\big(\gamma_1+\widetilde{L_1}(\ubm,\vbm,\xbm)\big)+L_1(\ubm,\vbm)+C_1(\ubm,\vbm,\xbm)\\+\nu_1S_2(\zbm)+E_1(\ubm,\vbm,\xbm,\zbm)=0,\end{multline*}
    where $\gamma_1$ is in $\RR\setminus\{0\}$, $\nu_1$ is any real number, $L_1 $ and $\widetilde{L_1}$ are two linear functions, $C_1\in\mfrak^2$ is a function of order two, and $E_1$ is in $\mfrak^3$.
    \item If the smooth equation is linear in any of the $\xbm$, we rename one of those variables to $v_K$, and the equation is
    \begin{multline*}\eta v_K+ (z_1+\dots+z_K)\big(\gamma_2+\widetilde{L_2}(\ubm,\vbm,\xbm)\big)+L_2(\ubm,\vbm,\xbm)+C_2(\ubm,\vbm,\xbm)\\+\nu_2S_2(\zbm)+E_2(\ubm,\vbm,\xbm,\zbm)=0,\end{multline*}
    for $\eta\in\RR\setminus\{0\}$, $\gamma_1,\nu_2\in\RR$, two linear functions $L_2,\widetilde{L_2}$, a function of order two $C_2\in\mfrak^2$, and a function $E_2\in\mfrak^3$;
\end{enumerate}
Then, the other equation, which has to yield a Morse singularity after eliminating the necessary variables, must be of the form:
\begin{enumerate}
    \item For the first case
    \begin{multline*}q_1(\ubm,\vbm,\xbm)+\alpha_1 v_{K-1}+(z_1+\cdots+z_K)M_1(\ubm,\vbm,\xbm)+\mu_1 S_2(z_1,\dots,z_{K})\\+E_3(\ubm,\vbm,\xbm,\zbm)=0,\end{multline*}
    for a function $M_1\in \mfrak$, a function $E_3\in\mfrak^3$ and a function $q_1$ in the corresponding variables that does not contain a linear term in $v_{K-1}$, since it is already taken into account in the term $\alpha_1 v_{K-1}$. Also, $\mu_1$ must be non-zero, because otherwise this multiple point space would have Milnor number one but that homology would not be alternating, contradicting the hypothesis that the germ has a good real perturbation.
    \item Similarly, for the second case\begin{multline*}q_2(\ubm,\vbm,\xbm)+\alpha_2 v_K+(z_1+\cdots+z_K)M_2(\ubm,\vbm,\xbm)+\mu_2 S_2(z_1,\dots,z_{K})\\+E_4(\ubm,\vbm,\xbm,\zbm)=0,\end{multline*}
    for a function $M_2\in \mfrak$, a function $E_4\in\mfrak^3$ and a function $q_2$ in the corresponding variables that does not contain a linear term in $v_K$, since it is already taken into account in the term $\alpha_2 v_K$. Again, for the same argument we mentioned for $\mu_1$, $\mu_2$ must be non-zero.
\end{enumerate}
On the one hand, on both of the previous cases, it is clear that $z_1,\ldots,z_k$ cannot appear as nonzero linear terms, since they would not yield a singular variety. 
On the other hand, the statement of the theorem only concerns terms in degree $K+2$, and $E_i$ would appear in terms of higher degree, and since the variety is an $A_1$ singularity and it is 2-determined, it will not be affected by the $E_i$ terms, therefore we may assume $E_1,E_2,E_3,E_4$ are zero.
Moreover, for the same reason, $C_i$ and $q_i$ can be assumed to be polynomials of degree (at most) two and $M_i$ a linear function. 
Finally, $q_i$ cannot contain any constants, since $D^{K+1}(g)$ is not empty; and it cannot contain linear terms in $\xbm$ or $\zbm$, since it is singular.

After all this, we can take an equivalence in the source and rename functions so that $\gamma_1=\eta=\mu_1=\mu_2=1$ and the $(K+2)$-jet has the form
\begin{equation}\label{eq:morralla1}
\begin{aligned}
    &\left(\ubm;\  \vbm;\  \xbm; \sum_{i=1}^{K-1}z^iu_i  + z^K q_1(\ubm,\vbm,\xbm) + \alpha_1 z^K v_{K-1} + z^{K+1}M_1(\ubm,\vbm,\xbm) + z^{K+2};\right.\\ 
    & \left.\sum_{i=1}^{K-1}z^iv_i + z^K\big( L_1(\ubm,\vbm)+ C_1(\ubm,\vbm,\xbm)\big) + z^{K+1}\big(1+\widetilde{L_1}(\ubm,\vbm,\xbm)\big) +\nu_1 z^{K+2} \right)
    \end{aligned}
\end{equation}

or the form
\begin{equation}\label{eq:morralla2}
\begin{aligned}
     &\left(\ubm;\  \vbm;\  \xbm; \sum_{i=1}^{K-1}z^iu_i  + z^K q_2(\ubm,\vbm,\xbm) + \alpha_2 z^K v_K + z^{K+1}M_2(\ubm,\vbm,\xbm) + z^{K+2};\right.\\ 
    & \left.\sum_{i=1}^{K-1}z^iv_i +z^K\big( v_K+ L_2(\ubm,\vbm,\xbm)+ C_2(\ubm,\vbm,\xbm)\big) +z^{K+1}\big(\gamma_2+\widetilde{L_2}(\ubm,\vbm,\xbm)\big)+\nu_2 z^{K+2} \right).
    \end{aligned}
\end{equation}

In the target, we denote the coordinates by $(\Ubm,\Vbm,\Xbm,Z_1,Z_2)\in\CC^{n+1}$ in the obvious way. 
We now show that we can take $\nu_1=\nu_2=0$.
To do this, taking the change of coordinates in the target $\psi(\Ubm,\Vbm,\Xbm,Z_1,Z_2)=(\Ubm,\Vbm,\Xbm,Z_1-\nu_jZ_2,Z_2)$ for $j=1$ in the first case and $j=2$ in the second. We can now rename the functions and make a change in source and target to maintain the sum in the form $\sum_{i=1}^{K-1}z^iv_i$. 
More precisely, after the change in the target we will have some terms in the $Z_1$ coordinate of both cases of the form
\[\sum_{i=1}^{K-1}z^iv_i-\nu_j\sum_{i=1}^{K-1}z^iu_i,\]
so we may take for each $i=1,\ldots,K-1$ the change of coordinates in the source taking $v_i$ into $v_i + \nu_ju_i$ and then fix the corresponding change in $\Vbm$ by a change of coordinates in the target (to take $V_i$ into $V_{i}-\nu_jV_{i})$.
Now we can rename the functions $q_j,L_j,\tilde L_j, C_j, M_j$ to preserve the previous form while they all preserve the same properties.
We will do this repeatedly, so we will omit the explicit description of these changes of coordinates whenever it is clear they can be done.
\newline

Let us focus on the first case now. We currently have a $(K+2)$-jet of the form of \cref{eq:morralla1} with $\nu_1=0$. Now, since $1+\widetilde{L_1}(\ubm,\vbm,\xbm)$ is a unit and $(\ubm,\vbm,\xbm)=(\Ubm,\Vbm,\Xbm)$, we may make the change of coordinates in the target that takes $Z_2$ into $Z_2(1+\widetilde{L_1}(\Ubm,\Vbm,\Xbm))^{-1}$. After taking coordinate changes on $(\ubm; \vbm; \xbm)$ and $(\Ubm;\Vbm;\Xbm)$ and renaming $L_1,C_1$ to take the unit into account, the $K+2$ jet from \cref{eq:morralla1} is equivalent to
\begin{equation*}
\begin{aligned}
    \Bigg( \ubm;\  \vbm;\  \xbm; & \sum_{i=1}^{K-1}z^iu_i  + z^K q_1(\ubm,\vbm,\xbm) + \alpha_1 z^K v_{K-1} + z^{K+1}M_1(\ubm,\vbm,\xbm) + z^{K+2};\\ 
    & \sum_{i=1}^{K-1}z^iv_i + z^K\big( L_1(\ubm,\vbm)+ C_1(\ubm,\vbm,\xbm)\big) + z^{K+1}  \Bigg).
    \end{aligned}
\end{equation*}

Now, if we perform the source change taking $z$ to
\[z-\binom{K+1}{1}^{-1}(L_1+C_1),\]
we eliminate the term $z^K(L_1+C_1)$ in $Z_2$ at the price of introducing some extra terms in the variables $\ubm,\vbm,\xbm$ to the coefficients of $z$ in degrees $1,\ldots,K-1$. 
These extra terms can be removed after an $\Ascr$-equivalence, although these changes also introduce new terms in $Z_1$ that can be included in $q_1,\alpha_1$ and $M_1$ without changing their properties. 
Finally, we can easily assume that $M_1=0$ by a target change, using the fact that we have the term $z^{K+1}$ in $Z_2$, and renaming $q_1$ and $\alpha_1$ if necessary, since the properties of $q_1$ will be preserved.
We are now left with the $(K+2)$-jet
\begin{equation*}
    \left( \ubm;\  \vbm;\  \xbm; \sum_{i=1}^{K-1}z^iu_i  + z^K q_1(\ubm,\vbm,\xbm) + \alpha_1 z^K v_{K-1} + z^{K+2};\sum_{i=1}^{K-1}z^iv_i + z^{K+1}  \right),
\end{equation*}
as desired. 
The properties claimed in the statement about $q_1=q$ and $\alpha_1=\alpha$ are proven with the same ideas we have shown in \cref{thm:K+2jet}: in order to have homology in the correct dimension so that the germ can have a good real perturbation, $q(\bm{0},\bm{0},\xbm)$ must be positive definite and $\alpha<1$ (recall also \cref{rem:althominreal}). Also, observe that a change of coordinates can take $q(\bm{0},\bm{0},\xbm)$ to the intended sum of squares.
\newline

We focus on the second case now, starting with \cref{eq:morralla2} where $\nu_2=0$. Let us remind what is the last pair of divided differences to compute $D^{K+1}(g)$:
\begin{equation}\label{eq:horroroso}\begin{cases}
q_2(\ubm,\vbm,\xbm)+\alpha_2 v_K+(z_1+\cdots+z_K)M_2+ S_2(z_1,\dots,z_{K})+E_4&=0,\\
v_K+ (z_1+\dots+z_K)(\gamma_2+\widetilde{L_2})+L_2+C_2+E_2&=0.
\end{cases}
\end{equation}

Observe that, if $\gamma_2\neq0$, we are returning to the previous case (or we get a contradiction). Indeed, if $\gamma_2$ is not zero, $\alpha_2$ must be zero in order for $D^{K+1}(g)$ to be singular as we are assuming, otherwise we would have a linear term $\alpha_2v_K$ on the first divided difference and another linear term $v_K+\gamma_2(z_1+\cdots+z_K)$. However, if $\alpha_2=0$, then we can eliminate one of the $z_i$ instead of $v_K$ and then $v_K$ must appear quadratically in $q_2(\bm{0},0,\dots,0,v_K,\xbm)$, since $v_K$ was never eliminated and that equation must give an $A_1$ singularity. This is, indeed, the same situation we had in the first case, simply naming $x_{n-2K+1}\coloneqq v_K$. Therefore, we assume that $\gamma_2=0$ and we have
\begin{equation*}
\begin{aligned}
    \Bigg( \ubm;\  \vbm;\  \xbm; &\sum_{i=1}^{K-1}z^iu_i  + z^K q_2(\ubm,\vbm,\xbm) + \alpha_2 z^K v_K + z^{K+1}M_2(\ubm,\vbm,\xbm) + z^{K+2};\\ 
    & \sum_{i=1}^{K-1}z^iv_i +z^K\big( v_K+ L_2(\ubm,\vbm,\xbm)+ C_2(\ubm,\vbm,\xbm)\big) +z^{K+1}\big(\widetilde{L_2}(\ubm,\vbm,\xbm)\big)\Bigg).
    \end{aligned}
\end{equation*}
Now we show that, possibly after a convenient source and target coordinate change, $\widetilde{L_2}$ must be a function such that $\widetilde{L_2}(\bm{0},\bm{0},\xbm)\neq0$. 
We have not yet used any information about $D^{K+2}(g)$, which must be empty, smooth or an $A_1$ of the correct dimension (this is satisfied as long as there are no divided differences that are identically zero, after eliminating variables). 
Indeed, $D^{K+2}(g)$ has a defining set of equations given by those of $D^{K+1}(g)$ stated before and the pair
\begin{equation}\label{eq:horror}\begin{cases}
\widetilde{L_2}(\ubm,\vbm,\xbm)+\widetilde{E_2}(\ubm,\vbm,\xbm,\zbm)&=0\\
M_2(\ubm,\vbm,\xbm)+S_1(z_1,\dots,z_{K+1})+\widetilde{E_4}(\ubm,\vbm,\xbm,\zbm)&=0,
\end{cases}\end{equation}
where $\widetilde{E_2},\widetilde{E_4}\in\mfrak^2$. Both equations have to contain linear terms different from $\ubm$ and $\vbm$, since they appeared before, in order to have an $A_1$ singularity. If   $\widetilde{L_2}$ is identically zero, that is not possible. If $\widetilde{L_2}(\bm{0},\bm{0},\xbm)=0$ but $\widetilde{L_2}$ is not identically zero, the only option left is that $\widetilde{L_2}$ has a linear term on $v_K$ and, at the same time, the function $L_2(\ubm,\vbm,\xbm)$ that appeared in the previous divided difference in \cref{eq:horroroso} has a linear term in $\xbm$ (say, $x_1$), so that eliminating $v_K$ in \cref{eq:horror} returns a linear term on $x_1$. However, then $D^{K+1}(g)$ would actually be smooth unless $\alpha_2=0$, for the same reason (see \cref{eq:horroroso}). So, we have two possibilities:
\begin{enumerate}
    \item[2.1] $\alpha_2=0$, $L_2$ has a linear term on $x_1$ and $\widetilde{L_2}$ is such that $\widetilde{L_2}(\bm{0},\bm{0},\xbm)=0$ but has a linear term on $v_K$; or
    \item[2.2] $\widetilde{L_2}(\bm{0},\bm{0},\xbm)\neq 0$, so it is a linear function with some $\xbm$ (which can be chosen to be $x_1$).
\end{enumerate}
In the first case, we may perform the change of coordinates in the source taking $v_K$ to be 
\[v_K-L_2(\ubm,\vbm,\xbm)-C_2(\ubm,\vbm,\xbm),\]
and, after that and a renaming of the functions $q_2,M_2$ and $\widetilde{L_2}$, we have a germ of the same form but with $L_2=C_2=0$ and such that we are in the second case $2.2$. This is because $\widetilde{L_2}$ had $v_K$ linearly and, since the change involves $L_2$ that had $x_1$ linearly, now the renamed $\widetilde{L_2}$ will have a linear term on $x_1$. Now, we perform the source change to remove everything but $x_1$ in $\widetilde{L_2}$, since it appears linearly and we can obtain (after a similar change involving $v_K$ if necessary in the case 2.2) a jet of the form
\begin{equation*}
\begin{aligned}
    \Bigg( \ubm;\  \vbm;\  \xbm;& \sum_{i=1}^{K-1}z^iu_i  + z^K q_2(\ubm,\vbm,\xbm) + \alpha_2 z^K v_K + z^{K+1}M_2(\ubm,\vbm,\xbm) + z^{K+2};\\ 
    & \sum_{i=1}^{K-1}z^iv_i +z^Kv_K +z^{K+1}x_1\Bigg).
    \end{aligned}
\end{equation*}

We can now eliminate $M_2$ by the source change that takes $z$ into
\[z-\binom{K+2}{1}^{-1}M_2(\ubm,\vbm,\xbm),\]
which eliminates $z^{k+1}$ and we can rename $q_2$ to
\[q_2+\binom{K+2}{2}\left(\binom{K+2}{1}^{-1}M_2(\ubm,\vbm,\xbm)\right)^2-\binom{K+1}{1}\binom{K+2}{1}^{-1}M_2(\ubm,\vbm,\xbm) ^2,\]
which has the same properties as $q_2$, i.e., without linear terms on $v_K$. We also had to adjust the variables $\ubm$ and $\vbm$. Finally, we can also eliminate the term $\alpha_2z^Kv_K$, if necessary, by performing a similar change of variables on $z$:
\[z-\binom{K+2}{1}\alpha_2v_K.\]
The remaining jet is
\begin{equation*}
\begin{aligned}
    \left( \ubm;\  \vbm;\  \xbm; \sum_{i=1}^{K-1}z^iu_i  + z^K q_2(\ubm,\vbm,\xbm)  + z^{K+2}; \sum_{i=1}^{K}z^iv_i +z^{K+1}x_1\right).
    \end{aligned}
\end{equation*}
The properties of $q\coloneqq q_2$ stated in the result follow from a similar argument to that we have shown in \cref{thm:K+2jet} and in the previous case. Finally, observe that a change of coordinates can take $q(\bm{0},\bm{0},\xbm)$ to the intended sum of squares (this is possible since we can do this modifying $x_1$ only up to a coefficient).
\end{proof}

\section{Examples}\label{sec:examples}

In this section, some examples are given on how to apply the previous theorems to classify map germs with good real perturbation.

\begin{example}
In \cite{Breva_Giménez2024} it is shown that the map germs $(\CC^3,0)\to (\CC^4,0)$ given by
\begin{align*}
    A_1 & \equiv (x,y,z(x^2+y^2+z^2),z^2), \\
    P_1 & \equiv (x,y,yz+z^4,xz+z^3), \\
    Q_2 & \equiv (x,y,z^3+y^2z,xz+yz^2), 
\end{align*}
are the only simple germs in the classification list of \cite{Houston1999a} which have a good real picture.

Notice that since $A_1$ and $Q_2$ are $3$-determined and $P_1$ is $4$-determined, together with \cref{thm:necessary}, they are actually the only possible map germs of corank one with a good real perturbation in these dimensions.
$A_1$ corresponds to the first case of \cref{thm:necessary} with $K=1$, $P_1$ is also the first case with $K=2$ and $Q_2$ is the second case of the same theorem for $K=1$.
\end{example}

\begin{example}
The corank one map germs $(\CC^4,0)\to(\CC^5,0)$ admitting good real perturbations must have one of the following jets (3,4 and 3-jets respectively):
\begin{align*}
    &(x_1,x_2,x_3,z(x_1^2+x_2^2+x_3^2) +z^3,z^2),\\
    &(u,v,x,zu+z^2(q(u,v,x)+\alpha v) + z^4, zv+z^3),\\
    &(v,x_1,x_2,z\tilde q(v,x_1,x_2) +z^3, vz+z^2x_1),
\end{align*}
where $\alpha <1$, and $q,\tilde q$ two functions without linear terms in $v$ and $q(0,0,x)=x^2$ and $q(0,x_1,x_2)=x_1^2+x_2^2$.
The first jet corresponds to the first case of \cref{thm:necessary} with $K=1$, the second jet corresponds to the same case with $K=2$, and the third corresponds to the second case of the same theorem with $K=1$.

\end{example}

\bigskip

{\bf Data availavility statement:} The authors declare that there is no associated data to this manuscript.

\bibliographystyle{myalpha.bst}
\bibliography{FullBib.bib}
\end{document}